\magnification\magstephalf
\input AHTOH-E.STY
\hfuzz2.7pt
\def\K{{\cal K}}
\def\H{{\bf H}}
\def\QQ{{\cal Q}}

\UDC{
512.543.72
+512.543.73
%
+512.544.33
+512.543.56  
}


\title{%
Equations over solvable groups
}
\author{%
Anton A. Klyachko$^{\flat\sharp}$
\quad
Mikhail A. Mikheenko$^{\flat\sharp}$
\quad
Vitaly A. Roman'kov$^{\natural\mathchar\times}$
}
\address{
$^\flat$%
Faculty of Mechanics and Mathematics of Moscow State University,
Moscow 119991, Leninskie gory, MSU.
\\
$^\sharp$Moscow Center for Fundamental and Applied Mathematics.
\\
$^\natural$Omsk branch of
Sobolev Institute of Mathematics,
Omsk 644043, ul. Pevtsova, 13.
\\
$^\mathchar\times$Dostoevsky Omsk State University,
Omsk 644077, prospekt Mira, 55-a.
\\
\vphantom{$\sum^\sum$}
klyachko@mech.math.msu.su
\qquad
mamikheenko@mail.ru
\qquad
romankov48@mail.ru
}
\grantsFirstSecond{
\RSF 22-11-00075}
\grantsThird{\RSF 22-21-00745}

\abstract{%
Not any nonsingular equation over a
metabelian
group has
solution in a larger
metabelian
group. However,
any nonsingular equation over
a solvable group
with
a subnormal
series with abelian \spacing{torsion-free} quotients
has a solution in a larger group with
a similar subnormal series of the same length
(and an analogous fact is valid for systems of equations).
}

\s 0.
Introduction

A system of equations
$\{w_i=1\;|\;i\in I\}$ with coefficients from a group $G$,
where $w_i$ are words in the alphabet~$G\sqcup X^{\pm1}$,
and $X$ is a set (the \emph{set of unknowns}),
is called
\emph{solvable over} $G$, if
there exists a group $\~G$
containing $G$ as a subgroup and a retraction of the free
product~$\~G*F(X)$ onto $\~G$ containing all elements~$w_i$ in its
kernel (henceforth, $F(X)$ is the free group with basis $X$). If
the \emph{solution group $\~G$} can be chosen from a class~$\K$,
then we say that the system \emph{is solvable in $\K$}.

The study of
solvability of equations over groups has a long history:
see, e.g.,
[GR62],
[Le62],
[B80],
[Ly80],
[How81],
[B84],
[EH91],
[How91],
[K93],
[KP95],
[FeR96],
[K97],
[K99],
[CG00],
[EdJu00],
[IK00],
[Juh\'a03],
[K06],
[P08]
[BK12],
[KL12],
[KT17],
[Ro17],
[BE18],
[ABA21],
[EH21],
[NT22],
[KM22],
and references therein;
see also surveys [Ro12], [NRR79], and book [LS80].

A system of equations (possibly infinite and with, possibly, infinitely
many unknowns) over a group is called \emph{nonsingular} if the
rows composed of the exponent-sums of unknowns in
each equation are linearly independent over~$\Q$. If these rows are 
linearly independent over the $p$-element field $\F_p=\Z/p\Z$ for each 
prime~$p$, then we call the system of equations \emph{unimodular}. In 
particular, one equation with one unknown is
\-
nonsingular if the exponent sum
of the unknown in this equation is nonzero;
\-
unimodular if this sum is $\pm1$.

\enditem
Unimodular equations behave better than arbitrary nonsingular ones:
e.g.,
in [K93] (see also [FeR96]), it was proven that
\disp{\sl
any unimodular equation over a torsion-free group is solvable over this
group;
}%
it is unknown whether a similar statement is valid for arbitrary
nonsingular equations.

For nilpotent groups, everything is
simple: Shmel'kin's theorem [Sh67] says
(in particular) that
\disp{\sl
any finite nonsingular system of equations over a nilpotent
torsion-free group $G$ has a
\(unique\) solution
in a nilpotent group
$\~G\supseteq G$ of the same nilpotency class;
namely, $\~G$ is the completion of $G$
\(if the system is
unimodular, then the
torsion-free
condition can be dropped, and the unique
solution exists in $G$ itself\).\fn{}%
}%
\footnote{}{\noindent$^{*)}$
Note that a theorem of Kuz'min [Ku74]
(see also [Ku06]) says that \sl a finitely generated
metabelian group $G$
is residually nilpotent if and only if
any unimodular equation
with one unknown
over $G$ has at most one solution in $G$
\rm(and the same is true for finitely generated
abelian-by-nilpotent
groups
[Ku78]).
}%
For solvable groups, the situation is more complicated:
Section 1 contains examples showing
(in particular) that
\disp{\sl
there exists a unimodular equation
with one unknown
over a metabelian
group (that can be chosen finite or, on the contrary,
torsion-free\), which is not solvable
in any larger metabelian group.
}%
(An example on this subject in [Ro17] is incorrect.)

In 1981, Howie suggested the following generalisation
of the well-known Kervaire--Laudenbach conjecture.

\proclaim Howie conjecture {\rm[How81]}.
Any nonsingular system of equations over any group is solvable over
this group.

It is unknown whether this is true or false, but there are many partial
results on this subject (see references above).
We
suggest to call a class $\K$ of groups a \emph{Howie class} if any
nonsingular system of equations over each group from $\K$ has
a solution in a larger group from~$\K$.

The quasivariety generated by a Howie class is a Howie class too
(see Section 2).
Therefore, we are mainly interested in
\emph{Howie quasivarieties}
(i.e. Howie classes which are quasivarieties).
Examples of such quasivarieties are
\-
the variety of all abelian groups,
and also the quasivariety of all abelian torsion-free groups
(an easy exercise);
\-
the quasivariety generated by all finite groups
([GR62] + Malcev theorem
on local residual finiteness of linear groups
[Ma40]).

\enditem
Examples from Section 1 mentioned above show that
the variety of all metabelian
groups is not a Howie class.
Nevertheless, we prove the
following theorem in Section 3.

\proclaim Main theorem.
The class consisting of all extensions of groups from a Howie quasivariety
by abelian torsion-free groups is a Howie quasivariety too.

In particular,
since abelian groups
form a Howie quasivariety,
an obvious induction leads us to the following fact
on solvable groups:
\disp{\sl\narrower\narrower\narrower
for each positive integer $n$,
the class of groups $G$
admitting a subnormal
series
$
G=G_1\triangleright\dots\triangleright G_{n+1}=\1
$
whose all quotients $G_i/G_{i+1}$ are abelian, and all, except
maybe the last one, are torsion-free is a Howie quasivariety,
i.e. any nonsingular system of equations over any group from this class
has a solution in a larger group from this class.
}%
(We can drop the words ``except maybe the last one"
and obtain another valid fact;
indeed, we can start with a simple observation:
{\sl the trivial quasivariety $\bigl\{\1\bigr\}$ is a Howie class},
and then use a similar induction argument.)
Note, however, that the minimum length of such a series can be larger than
the derived length;
see Section~1.

The authors thank A. I. Budkin
and A. V. Khudyakov
for valuable remarks.
We are grateful to an anonymous referee whose comments 
allow us to improve the text.
The first two authors thank also the Theoretical Physics and Mathematics
Advancement Foundation ``BASIS".

\s 1.
Discouraging examples

\Lemma 1.
If elements $x$ and $y$ of a metabelian group are such that
$y^6$ and $u=xx^{-y}x^{y^2}$ belong to the commutator subgroup, then
$
uu^yu^{-y^3}u^{-y^4}=1.
$
\newline
\rm{(Henceforth,
we use standard abbreviations: $g^h\:=h^{-1}gh$,\quad
$g^{-h}\:=h^{-1}g^{-1}h$,\quad
$g^{\the\year h}\:=h^{-1}g^{\the\year}h$,
\dots)}

\Proof
The elements $u=xx^{-y}x^{y^2}$ and $x$
are equal modulo the commutator subgroup (in any group); hence,
$x$ lies in the commutator subgroup.
Thus,
$$
uu^yu^{-y^3}u^{-y^4}=
\(xx^{-y}x^{y^2}\)\(xx^{-y}x^{y^2}\)^y\(xx^{-y}x^{y^2}\)^{-y^3}\(xx^{-y}x^{y^2}\)^{-y^4}
=
xx^{-y^6}
=
1.
$$
Here, we use two facts:
\-
conjugates of $x$ commute
(because $x$ lies in the commutator subgroup of a metabelian group);
\-
$(1-y+y^2)(1+y-y^3-y^4)=1-y^6$ (which is an elementary identity).
\enditem
This completes the proof.

\Proposition 1.
There exists a metabelian group $G$ and a unimodular equation
$w(x)=1$
over $G$ having no solutions in any metabelian group
$\~G\supseteq G$.
Moreover, the group $G$ can be chosen
\item{\rm a)}
finite \(of order 42\/\)%
\fn{%
42 is the minimum possible order of such an example [M23].
}
\item{\rm b)}
or, on the contrary, torsion-free and such that any
nonsingular system of equations over $G$ has a solution
in a larger three-step solvable group.

\Proof
Let $G$ be a metabelian group with elements $a$ and $c$ such that
$a^6=1$ and $c$ is contained in the commutator subgroup of $G$.
By Lemma 1, the unimodular equation
$xx^{-a}x^{a^2}=c$
cannot
have solutions in larger metabelian
groups if $cc^ac^{-a^3}c^{-a^4}\ne1$.
A natural example of such a group is the
group of triangular
invertible $2\times2$-matrices over the group
ring~$\Z[y]/\bigl((y^6-1)\cdot\Z[y]\bigr)$
of the cyclic group of order six:
$$
a=
\pmatrix{
1&0
\cr
0&y
},\
c=
\pmatrix{
1&y-1
\cr
0&1
}
=
\[
\pmatrix{
1&1
\cr
0&1
}
,a
\]
\hbox{ and }
cc^ac^{-a^3}c^{-a^4}=
\pmatrix{
1&(y-1)(1+y-y^3-y^4)
\cr
0&1
}
\ne1.
$$
(The commutator is understood as $[u,v]\:=u^{-1}v^{-1}uv$.)

To construct a group of order 42, it suffices to make the same
construction
over the seven-element field~$\F_7=\Z/7\Z$, which is
a homomorphic image of the ring $\Z[y]/\bigl((y^6-1)\cdot\Z[y]\bigr)$,
$y\mapsto 5$
(five is the generator of the multiplicative group of
this field):
$$
\eqalign{
G=&\left\{
\pmatrix{
1&u
\cr
0&v
}
\;\Bigm|\;
u,v\in\F_7,\;v\ne0
\right\},
\cr
a=
\pmatrix{
1&0
\cr
0&5
},\
c=
\pmatrix{
1&5-1
\cr
0&1
}
=&
\[
\pmatrix{
1&1
\cr
0&1
}
,a
\]
\hbox{ and }
cc^ac^{-a^3}c^{-a^4}=
\pmatrix{
1&(5-1)(1+5-5^3-5^4)
\cr
0&1
}
\ne1.
}
$$


\medskip

To construct a torsion-free group, we modify the first
example:
$$
\eqalign{
G_1=&\left\{
2^k\pmatrix{
1&u
\cr
0&y^k
}
\;\Bigm|\;
u\in\Z[y]/\bigl((y^6-1)\cdot\Z[y]\bigr),\;k\in\Z
\right\},
\cr
a=
2\pmatrix{
1&0
\cr
0&y
},\
c=
\pmatrix{
1&y-1
\cr
0&1
}
=&
\[
\pmatrix{
1&1
\cr
0&1
}
,a
\]
\hbox{ and }
cc^ac^{-a^3}c^{-a^4}=
\pmatrix{
1&(y-1)(1+y-y^3-y^4)
\cr
0&1
}
\ne1.
}
$$
Clearly, the group $G_1$ is torsion-free, and $a^6$ lies in its centre.
The central product
$$
G=G_1
\bigtimes\limits_{a^6=f}
H,
\qqbox{where}
H=\pres<d,e,f|
[d,e]=f,\;[d,f]=[e,f]=1>,
$$
of $G_1$ and the Heisenberg group $H\iso\UT_3(\Z)$
satisfies conditions of Lemma 1 and is torsion-free
(because, in the Heisenberg group, the centre $\gp f$ is
\emph{isolated}, i.~e., if $g^k\in\gp f$ for some $k\in\Z\setminus\0$, then
$g\in\gp f$).
Moreover, $G$ has a subnormal series of length three
$$
G\triangleright G_1
\triangleright
\left\{
\pmatrix{
1&u
\cr
0&1
}
\;\Bigm|\;
u\in\Z[\gp y_6]
\right\}
\triangleright
\1
$$
whose quotients, obviously, are free abelian groups
of ranks two, one, and six. This completes the proof of Proposition 1,
because, by virtue of the main theorem, a group with a
length-three subnormal series with abelian torsion-free quotients form
a Howie class.

(We note parenthetically
that, in the last example, there are no length-two subnormal series with
abelian torsion-free quotients, because, by the main theorem, such a
series would imply the solvability of the equation $xx^{-a}x^{a^2}=c$ in a
larger metabelian group, which is not the case, as we have shown.)

\s 2.
Quasivarieties and systems of equations

Recall that a \emph{quasivariety} of groups is the class of all groups
satisfying a given (finite or infinite) set of
\emph{quasi-identities}, i.e. sentences of the form
$$
\forall x,y,\dots\
v_1(x,y,\dots)=1\;\&\dots\&\;v_k(x,y,\dots)=1
\imp
w(x,y,\dots)=1,
$$
where $v_i$ and $w$ are some elements of the free group $F(x,y,\dots)$.
For example, the quasivariety of abelian torsion-free groups
can be defined by the following quasi-identities
(the last one is, in fact, an identity):
$
\forall x\;
x^2=1\imp x=1,\quad
\forall x\;
x^3=1\imp x=1,\quad
\dots,\quad
\forall x,y\;
1=1\imp [x,y]=1.
$
Another (equivalent) definition of the notion of quasivariety
can be found below, in the beginning of the proof of Proposition 2.
For more details about quasivarieties, see the books
[Bu02] and [Go99]; we
recall only that
\disp{\sl\hfuzz5pt
the class consisting of all possible extensions
of groups from a quasivariety by groups from another
quasivariety is a quasivariety {\rm[Ma67]}
\rm(called the \emph{product} of two quasivarieties\).
}%
Therefore, to prove the main theorem, it suffices to show that the class
of extensions (in that theorem)
is a Howie class.

\Proposition 2.
For any quasivariety $\QQ$,
the class of groups
$$
\H({\QQ})
=
\left\{G\in\QQ\;\Biggm|\;
\vcenter{\rm
\hbox{any nonsingular system of equations over $G$}
\hbox{has a solution in a larger group from $\QQ$}
}
\right\}
$$
coincides with the class
$$
\H_{\rm fin}({\QQ})
=
\left\{G\in\QQ\;\Biggm|\;
\vcenter{\rm
\hbox{any finite nonsingular system of equations over $G$}
\hbox{has a solution in some larger group from $\QQ$}
}
\right\}
$$
and is a quasivariety.

\Proof
By a theorem of Malcev (see [Ma70], Theorem 5.11.4),
\disp{\sl
quasivarieties of groups are the same as nonempty classes of groups
closed with respect to subgroups and reduced (= filtered) products.
}%
Recall that the \emph{reduced} (or \emph{filtered}) product
$\(\!\!\bigtimes G_i\)/F$
of groups $G_i$, where $i\in I$, with respect to a filter~$F$
on the set~$I$ is the quotient
group
$
\(\bigtimes\limits_{i\in I}G_i\)\biggm/
\left\{h\in\bigtimes\limits_{i\in I}G_i\;\Biggm|\;
\{i\in I\;|\;h_i=1\}\in F\right\}.
$

The closedness of classes $\H({\QQ})$ and $\H_{\rm fin}({\QQ})$
with respect to these operations is fairly obvious
(to obtain a solution group for a system of equations over
a reduced product, we can take
the reduced product (with respect to the same filter)
of solution groups for the corresponding systems of equations over the
factors).

The coincidence of the classes $\H({\QQ})$ and $\H_{\rm fin}({\QQ})$
is a special case of the G\"odel--Malcev
compactness theorem; but for reader's convenience, we
give a sketch of a direct proof:

\-
    consider an (infinite) nonsingular system $V$ of equations
    over a group $G\in\H_{\rm fin}({\QQ})$ with
    an (infinite) set of unknowns $X$;

\-
    for each finite subsystem $W\subseteq V$, take its solution group
    $G_W\in\QQ$;

\-
    take the Cartesian product of all $G_W$;

\-
    and choose a set (of values of the unknowns) $\~X$ in this product:
    their coordinates $\~X_W$ form a solution to the finite system
    $W$;

\-
    on the set of finite subsystems $\{W\}$, consider the filter $F$
    consisting of all sets containing all
    ``sufficiently large systems":
    $
    F=\bigl\{M\;|\;
    \exists W\  M\supseteq\{W'\;|\;W'\supseteq W\}\bigr\};
    $
\-
    the group $G$ is embedded diagonally into this reduced product
    $\~G=\(\bigtimes\limits_W G_W\)/F\in\QQ$;

\-
    the set $\~X\subset\~G\in\QQ$ forms a solution to the entire system
    $V$,
    (because, for each equation~$w(X)=1$ from the system $V$, the
    coordinates $\~X_W$ form a solution to this equation if $W\ni w$
    (i.~e., $w(\~X_W)=1$); and this set of
    coordinates $\{W\;|\;W\ni w\}$ lies in the filter $F$);
    \-
    thus,
    $G\in\H({\QQ})$, and this completes the proof of Proposition 2.

\medskip
\noindent
Proposition 2 shows, in particular,
that
\item{1)}
the quasivariety generated by a Howie class
is a Howie class itself;
\newline
indeed,
for the quasivariety $\qvar\K$
generated by a Howie class $\K$, we obtain that
$\H(\qvar\K)$
{\itemitem{--}
is a quasivariety (by Proposition 2)
\itemitem{--}
and contains $\K$
(since $\K$ is a Howie class, i.e. nonsingular systems of equations
over groups from $\K$ are solvable in $\K\subseteq\qvar\K$);
\itemitem{--}
hence, $\H(\qvar\K)\supseteq\qvar\K$ (as
$\qvar\K$ is the least quasivariety containing $\K$);
\itemitem{--}
the last inclusion is an equality actually, because
the inverse inclusion holds by the definition of~$\H(\dots)$

}
(therefore, it is natural to study Howie quasivarieties
rather than
arbitrary Howie classes);
\item{2)}
in the proof of the main theorem, we can restrict ourselves to finite
systems of equations
(since the main theorem asserts the equality
$\H(\K)=\K$ for a quasivariety $\K$,
while Proposition 2 says that $\H(\K)=\H_{\rm fin}(\K)$),
i.~e., the main theorem is equivalent to the following
(formally weaker) statement.

\enditem
\proclaim Main theorem {\rm (a simplified form)}.
Suppose that a group $G$ contains a normal subgroup~$B$
from a Howie quasivariety
$\K$, and the quotient group $G/B=A$
is abelian and torsion-free. Then any finite
nonsingular system of equations
$
\{w_i=1\;|\; i\in I\}
$
{\(where $w_i\in G*F(X)$ and $I$ is a finite set\)}
has a solution in some group $\~G\supseteq G$,
which is also an extension of a group from $\K$ by an abelian
torsion-free group.

\s 3.
Proof of the main theorem

Over the group $A$, we have the induced system of equations:
$
\{\=w_i=1\;|\; i\in I\},
$
where $\=w_i\in A*F(X)$ are images
of $w_i$ under the natural homomorphism $G*F(X)\to A*F(X)$.

The first simple observation is that
\disp{\sl
we can assume that the induced system of equations
$
\{\=w_i=1\}
$
has a solution in $A$.}%
Indeed,
this system of equations over $A=G/B$
has solution in some larger torsion-free abelian group~$\~A$
(because the class of abelian torsion-free groups is a Howie class);
therefore, we can embed
$G$ into the Cartesian (= unrestricted) wreath product $B\bar{\wr}A$
by the Kaloujnine--Krasner theorem [KK51],
and then embed the wreath product $B\bar{\wr}A$
into $B\bar{\wr}\~A$; now, we replace
$A$ with $\~A$, and $B$ with the base of the wreath product
(using closedness of the quasivariety $\K$ with respect to Cartesian
powers); we obtain what we want.
As a byproduct of this
argument, we come to the conclusion that
\disp{\sl
we can assume that the initial extension
splits:  $G=A\semitimes B$.
}%

The next observation is that
\disp{\sl\narrower\narrower\narrower\narrower
we can assume that a solution to the induced system of
equations
$
\{\=w_i=1\}
$
in $A$ is the set of identity elements.}%
Indeed, this can be achieved by an obvious change of variables
(because we assume already that some solution exists in $A$).

This observation means that $w_i$ can be rewritten
as words in the alphabet
$B\sqcup\bigsqcup\limits_{a\in A}X^{\pm a}$:
$$
\eqalign{
&G*F(X)\ni w_i=f(v_i),
\qbox{where $v_i\in B*F\(\bigsqcup\limits_{a\in A}X_a\)$ are
some words, $X_a$ --- copies of the alphabet $X$},
\cr
&\qbox{and the
homomorphism
$f\:B*F\(\bigsqcup\limits_{a\in A}X_a\)\to G*F(X)$
sends $X_a\ni x_a$ to $x^a\:=a^{-1}xa$ and $B\ni b$ to $b$.}
}
$$
On the group $B*F\(\bigsqcup\limits_{a\in A}X_a\)$,
the group $A$ acts naturally (on the right):
$b\o a=a^{-1}ba$ and $x_{a'}\o a=x_{a'a}$.
The key observation
(whose proof is postponed,
see Lemma 2 below) is that
\dispno{\sl
the system of equations $\{v_i\o a=1\;|\;i\in I,\;a\in A\}$
over $B$
\(with unknowns $\bigsqcup\limits_{a\in A}X_a$\)
is nonsingular.
}(*)%
Now everything is simple.
\-
Since the system of equations $\{v_i\o a=1\;|\;i\in I,\;a\in A\}$
over $B$
is nonsingular, it has a solution in a group $\~B\supseteq B$
from $\K$ (because $\K$ is a Howie class);
\-
i.e. the quotient group
$
U=
\(B*F\(\bigsqcup\limits_{a\in A}X_a\)\)\Big/
\nc{\{v_i\o a=1\;|\;i\in I,\;a\in A\}}
$
admits a homomorphism to $\~B$ whose restriction to $B$ is the
identity map.
\-
Therefore, the quotient group $\^B\supseteq B$ of $U$
by the intersection of the kernels of all
homomorphisms to~$\~B$
\itemitem{--}
contains $B$ and a solution to the system of
equations
$\nc{\{v_i\o a=1\;|\;i\in I,\;a\in A\}}$ over $B$
(obviously);
\itemitem{--}
is contained in $\K$
(because $\^B$ is sub-Cartesian product of subgroups of
$\~B\in\K$, and the class $\K$
is a
quasivariety).
\-
It remains to note that the natural action of $A$ on $U$
(which is well defined, since
the subgroup $\nc{\{v_i\o a=1\;|\;i\in I,\;a\in A\}}$ is $A$-invariant)
extending the action of $A$ on $B$
induces an action of $A$ on $\^B$
(because the intersection of the
kernels of all homomorphisms from a group to another one
is a characteristic subgroup).
\-
Therefore, the corresponding semidirect product
$\~G=A\semitimes\^B\supseteq A\semitimes B=G$ is the required group
containing a solution to the initial system and being an extension of
a group from $\K$ by an abelian torsion-free group.

\enditem
It remains to prove $(*)$.
This fact is a corollary of the following lemma for
\-
$R=\Z[A]$ is the group ring of $A$
(this ring has no zero divisors, as is well known ---
we leave a proof to the reader as an easy exercise),
\-
$M$ is the free abelian group with
basis $\bigsqcup\limits_{a\in A}X_a$
or, which is the same thing, $M$ is the free $R$-module with basis $X$
(where $ax=x_a$),

\-
$\{m_1,\dots,m_k\}\subset M$ is the image of an arbitrary finite
subset of $\{v_i\;|\;i\in I\}$ under the natural homomorphism
$B*F\(\bigsqcup\limits_{a\in A}X_a\)\to M$
with the kernel $(\hbox{the commutator subgroup})\cdot B$
(in other words, $m_i$ are tuples of exponent-sums of
the unknowns 
(from $\bigsqcup\limits_{a\in A}X_a$)	
in the equations $v_i=1$),
\-
$I$ is the \emph{fundamental ideal},
i.e. the kernel of the homomorphism $\Z[A]\to\Z$
sending $A$ to 1
(i.~e., $M/(IM)$, as a module over $R/I\iso\Z$,  
is naturally isomorphic to the free $\Z$-module
with basis $X$, the elements $m_i+IM\in M/(IM)$ correspond
to the
tuples of exponent-sums of the unknowns from  
$X$ in the equations $w_i=1$).

\Lemma 2.
Suppose that elements $m_1,\dots,m_k$ of a free module $M$ over
an associative commutative ring $R$
without zero divisors are
linearly dependent over $R$.
Then, for any prime ideal $I\nin R$, the
elements~$m_i+IM$ of the quotient module $M/(IM)$
over the quotient ring $R/I$
are linearly dependent
over $R/I$.

\Proof
Let us choose some basis $e_1,\dots,e_n$ of the module $M$ (which
can surely be assumed finitely generated,
because any finite subset of a free module is contained in
a finitely generated free submodule).
The linear dependence
of elements $m_i=\sum\limits_j r_{ij}e_j$ means that all the minors
of order $k$ in the matrix $(r_{ij})$ are zero (because the ring $R$
embeds into a field,
and, over a field, rows of length $k$ are linearly dependent
if and only if all $k$-minors vanish, as is well known).
Therefore, all minors of order $k$ in the matrix
$(r_{ij}+I)$ are zero. Therefore,
$m_i+IM=\sum\limits_j(r_{ij}+I)(e_j+IM)$
are linearly dependent (because the ring $R/I$ embeds into a field too).
This completes the proofs of the lemma and the main theorem.

\baselineskip 11.8pt
\References

[ABA21]
M. F. Anwar, M. Bibi, M. S. Akram,
On solvability of certain equations of arbitrary length
over torsion-free groups,
Glasgow Mathematical Journal, 63:3 (2021), 651-659.
\arXiv 1903.06503

[BK12]
D. V. Baranov, A. A. Klyachko,
Economical adjunction of square roots to groups,
Siberian Math. Journal, 53:2 (2012), 201-206.
\arXiv 1101.3019

[BE18]
M. Bibi, M. Edjvet,
Solving equations of length seven over torsion-free groups,
Journal of Group Theory, 21:1 (2018), 147-164.

[B80]
S. D. Brodskii,
Equations over groups and group with a single defining relation,
Russian Math. Surveys, 35:4 (1980), 165-165.


[B84]
S. D. Brodskii,
Equations over groups and group with one relator,
{Siberian Math. Journal}, 25:2 (1984), 235-251.


[Bu02]
A. I. Budkin,
Quasivarieties of groups [in Russian],
Altay State Univ., Barnaul, 2002.

[CG00]
A. Clifford, R. Z. Goldstein,
Equations with torsion-free coefficients,
{Proc. Edinburgh Math. Soc.}, {43:2} (2000), 295-307.

[EH91]
M. Edjvet, J. Howie,
The solution of length four equations over groups,
{Trans. Amer. Math.  Soc.}, {326:1} (1991), 345-369.

[EH21]
M. Edjvet, J. Howie,
On singular equations over torsion-free groups,
International Journal of Algebra and Computation, 31:3 (2021), 551-580.
\arXiv:2001.07634

[EdJu00]
M. Edjvet, A. Juh{\accent 19 a}sz,
Equations of length 4 and one-relator products,
{Math. Proc. Cambridge Phil. Soc.}, 129:2 (2000), 217-230.

[FeR96]
R. Fenn, C. Rourke,
Klyachko's methods and the solution of
equations over torsion-free groups,
\newline
{L'Enseignment
Math{\accent 19 e}matique}, 42 (1996), 49-74.

[GR62]
M. Gerstenhaber, O.S. Rothaus,
The solution of sets of equations in groups,
{Proc. Nat. Acad. Sci. USA}, {48:9} (1962), 1531-1533.



[Go99]
V.A. Gorbunov,
Algebraic theory of quasivarieties.
Siberian School of Algebra and Logic.
Consultants Bureau,
New York, 1998.

[How81]
J. Howie,
On pairs of 2-complexes and systems of equations over groups,
{J. Reine Angew Math.}, 1981:324 (1981), 165-174.

[How91]
J. Howie,
The quotient of a free product of groups by a single high-powered relator.
III: The word problem,
{Proc. Lond. Math. Soc.}, {62}:3 (1991), 590-606.

[IK00]
S. V. Ivanov, A. A. Klyachko,
Solving equations of length at most six over torsion-free groups,
Journal of Group Theory, 3:3 (2000), 329-337.

[Juh{\accent 19 a}03]
Juh{\accent 19 a}sz A.
On the solvability of a class
of equations over groups,
{Math. Proc. Cambridge Phil. Soc.}, 135:2 (2003), 211-217.

[K93]
A. A. Klyachko,
A funny property of sphere and equations over groups,
Communications in Algebra, 21:7 (1993), 2555-2575.

[K97]
A. A. Klyachko,
Asphericity tests,
International Journal of Algebra and Computation, 7:4 (1997), 415-431.

[K99]
A. A. Klyachko,
Equations over groups, quasivarieties,
and a residual property of a free group,
Journal of Group Theory, 2:3 (1999), 319-327.

[K06]
A. A. Klyachko,
How to generalize known results on equations over groups,
Math. Notes, 79:3 (2006), 377-386.
\arXiv math.GR/0406382

[KL12]
A. A. Klyachko, D. E. Lurye,
Relative hyperbolicity and similar properties of one-generator one-relator
relative presentations with powered unimodular relator,
Journal of Pure and Applied Algebra, 216:3 (2012), 524-534.
\arXiv 1010.4220

[KM22]
A. A. Klyachko, M. A. Mikheenko,
Yet another Freiheitssatz: Mating finite groups with locally indicable ones,
Glasgow Mathematical Journal (to appear).
\arXiv 2204.01122

[KP95]
A. A. Klyachko, M. I. Prishchepov,
The descent method for equations over groups,
Moscow Univ. Math. Bull. 50:4 (1995), 56-58.

[KT17]
A. A. Klyachko, A. B. Thom,
New topological methods to solve equations over groups,
Algebraic and Geometric Topology, 17:1 (2017), 331-353.
\arXiv 1509.01376

[KK51]
M. Krasner, L. Kaloujnine,
Produit complet des groupes de permutations et le
probl\`eme d'extension de groupes. III,
Acta Sci. Math., 14 (1951), 69-82.

[Ku74]
Yu. V. Kuz'min,
Residual properties of metabelian groups [in Russian],
Algebra i Logika, 13:3 (1974), 300-310.

[Ku78]
Yu. V. Kuz'min,
Some approximation properties the variety $\eufm{AN}_c$,
Russian Math. Surveys, 33:4 (1978), 257-258.

[Ku06]
Yu. V. Kuz'min,
Homological group theory [in Russian].
Moscow: Faktorial, 2006.

[Le62]
F. Levin,
Solutions of equations over groups,
{Bull. Amer. Math. Soc.}, {68:6} (1962), 603-604.

[Ly80]
R. C. Lyndon,
Equations in groups,
{Bol. Soc. Bras. Math.}, {11}:1 (1980), 79-102.

[LS80]
R. Lyndon, P. Schupp,
Combinatorial group theory.
Springer, 1977.

[Ma40]
A. Malcev,
On isomorphic matrix representation of infinite groups [in Russian],
Mat. Sb., 8(50):3 (1940), 405-422.

[Ma67]
A. I. Maltsev,
Multiplication of classes of algebraic systems,
Siberian Math. J., 8:2 (1967), 254-267.

[Ma70]
A. I. Mal'cev,
Algebraic systems.
Springer-Verlag, New York, Heidelberg, and Berlin, 1973.

[M23]
M. A. Mikheenko,
On $p$-nonsingular systems of equations over solvable groups,
arXiv:2309.09096\thinspace.

[NT22]
M. Nitsche, A. Thom,
Universal solvability of group equations,
Journal of Group Theory, 25:1 (2022), 1-10.
\arXiv 1811.07737

[NRR79]
G. A. Noskov, V. N. Remeslennikov, V. A. Roman'kov,
Infinite groups,
J. Soviet Math., 18:5 (1982), 669-735.

[P08]
V. G. Pestov,
Hyperlinear and sofic groups: A brief guide,
Bull. Symb. Log., 14:4 (2008), 449-480.
\arXiv:0804.3968

[Ro12]
V. A. Roman'kov,
Equations over groups,
Groups-Complexity-Cryptology, 4:2 (2012), 191-239.

[Ro17]
V. A. Roman'kov,
On solvability of regular equations in the variety of metabelian
groups,
Applied discrete mathematics, 36 (2017), 51-58.

[Sh67]
A. L. Shmel'kin,
Complete nilpotent groups [in Russian],
Algebra i Logika Seminar, 6:2 (1967), 111-114.

\end